\documentclass[final]{siamltex}

\usepackage[hypertex]{hyperref}
\usepackage{amsmath,amssymb}

\title{Singular trajectories of control-affine systems}

\author{Yacine Chitour\footnotemark[2]
\and Fr\'ed\'eric Jean\footnotemark[3]
\and Emmanuel Tr\'elat\footnotemark[4]}

\newtheorem{remark}[theorem]{Remark}

\renewcommand{\leq}{\leqslant}
\renewcommand{\geq}{\geqslant}

\newcommand{\beq}{\begin{equation}}
\newcommand{\eeq}{\end{equation}}
\newcommand{\bea}{\begin{eqnarray}}
\newcommand{\eea}{\end{eqnarray}}
\newcommand{\beaa}{\begin{array}}
\newcommand{\eeaa}{\end{array}}

\newcommand{\bed}{\begin{itemize}}
\newcommand{\eed}{\end{itemize}}
\newcommand{\bedd}{\begin{description}}
\newcommand{\eedd}{\end{description}}

\def\N{\mathbb{N}}
\def\R{\mathbb{R}}

\begin{document}
\maketitle

\renewcommand{\thefootnote}{\fnsymbol{footnote}}

\footnotetext[2]{Labo.\ des Signaux et Syst\`emes, Universit\'e
  Paris-Sud, CNRS, Sup\'elec, 91192 Gif-sur-Yvette cedex, France ({\tt
    Yacine.Chitour@lss.supelec.fr}).}
\footnotetext[3]{ENSTA, UMA, 32 bld Victor, 75739 Paris, France (\tt
  Frederic.Jean@ensta.fr).}
\footnotetext[4]{Universit\'e Paris-Sud, Math., UMR 8628, Bat.\ 425,
  91405 Orsay cedex, France ({\tt Emmanuel.Trelat@math.u-psud.fr}).}

\renewcommand{\thefootnote}{\arabic{footnote}}

\begin{abstract}
  When applying methods of optimal control to motion planning or
  stabilization problems, some theoretical or numerical difficulties
  may arise, due to the presence of specific trajectories, namely,
  singular minimizing trajectories of the underlying optimal control
  problem. In this article, we provide characterizations for singular
  trajectories of control-affine systems. We prove that, under generic
  assumptions, such trajectories share nice properties, related to
  computational aspects; more precisely, we show that, for a generic
  system -- with respect to the Whitney topology --, all nontrivial
  singular trajectories are of minimal order and of corank one. These
  results, established both for driftless and for control-affine
  systems, extend results of \cite{CJTCRAS,CJT}.  As a consequence,
  for generic systems having more than two vector fields, and for a
  fixed cost, there do not exist minimizing singular trajectories. We
  also prove that, given a control system satisfying the LARC,
  singular trajectories are strictly abnormal, generically with
  respect to the cost. We then show how these results can be used to
  derive regularity results for the value function and in the theory
  of Hamilton-Jacobi equations, which in turn have applications for
  stabilization and motion planning, both from the theoretical and
  implementation issues.
\end{abstract}

\section{Introduction}
When addressing standard issues of control theory such as motion
planning and stabilization, one may adopt an approach based on optimal
control, e.g., Hamilton-Jacobi type methods and shooting
algorithms. One is then immediately facing intrinsic difficulties due
to the possible presence of singular trajectories. It is therefore
important to characterize these trajectories, by studying in
particular their existence, optimality status, and the related
computational aspects.
In this paper, we provide answers to the aforementioned questions for
control-affine systems, under generic assumptions, and then
investigate consequences in optimal control and its applications.

Let $M$ be a smooth (i.e.\ $C^{\infty}$) manifold of dimension $n$.
Consider the control-affine system
$$
(\Sigma)\qquad\dot{x}=f_0(x)+\sum_{i=1}^m u_if_i(x),
$$
where $x\in M$, $m$ is a positive integer, $(f_0,\ldots,f_m)$ is a 
$(m+1)$-tuple of smooth vector fields on $M$, and the control 
$u=(u_1,\ldots,u_m)$ takes values in an open subset $\Omega$ of
$\mathbb{R}^m$.
For $x_0\in M$ and $T>0$, a control $u\in L^\infty([0,T],\Omega)$ is
said to be \textit{admissible} if the trajectory $x(\cdot,x_0,u)$ of
$(\Sigma)$ associated to $u$ and starting at $x_0$ is well defined on
$[0,T]$. On the set ${\cal U}_{x_0,T}$ of
admissible controls, define the {\em end-point
  mapping} by
$$E_{x_0,T}(u) := x(T,x_0,u).$$
It is classical that ${\cal U}_{x_0,T}$ is an open subset of
$L^\infty([0,T],\Omega)$ and that $E_{x_0,T}:{\cal U}_{x_0,T} \rightarrow
M$ is a smooth map.

\begin{definition}
A control $u\in {\cal U}_{x_0,T}$ is said to be {\em singular} if $u$
is a critical point of the end-point mapping $E_{x_0,T}$,
i.e.\ its differential at $u$, $DE_{x_0,T}(u)$, is not surjective.
A trajectory $x(\cdot,x_0,u)$ is said to be {\em singular} if
$u$ is singular and {\em of corank one} if the codimension in the tangent
space of the range of $E_{x_0,T}(u)$ is equal to one.
\end{definition}

In other words, a control $u\in {\cal U}_{x_0,T}$ is singular if the
linearized system along the trajectory $x(\cdot,x_0,u)$ is not
controllable on $[0,T]$.
Singular trajectories appear as singularities in the set of
trajectories of $(\Sigma)$ joining two given points, and hence,
they play a crucial role in variational problems associated to
$(\Sigma)$ and in optimal control, as described next.

Let $x_0$ and $x_1$ be two points of $M$, and $T>0$.
Consider the following optimal control problem:
among all the trajectories of $(\Sigma)$ steering $x_0$ to $x_1$,
determine a trajectory minimizing the {\em cost}
\begin{equation}\label{coutintro}
C_{U,\alpha,g}(T,u)=\int_0^T\Big(\frac{1}{2}
u(t)^TU(x(t))u(t)+\alpha(x(t))^Tu(t) 
+ g(t,x(t)) \Big) dt,
\end{equation}
where $\alpha=(\alpha_1,\ldots,\alpha_m)\in C^{\infty}(M,\R^m)$,
$g\in C^\infty(\R\times M)$, and
$U$ takes values in the set of symmetric positive definite $m\times
m$ matrices.

According to the Pontryagin Maximum Principle (see \cite{P}), for every
optimal trajectory $x(\cdot):=x(\cdot,x_0,u)$, there exists a nonzero
pair $(\lambda(\cdot),\lambda^0)$, where $\lambda^0$ is
a nonpositive real number and $\lambda(\cdot)$ is an
 absolutely continuous function on $[0,T]$ (called {\em
adjoint vector}) with $\lambda(t) \in T^*_{x(t)}M$,
such that, almost everywhere on $[0,T]$,
\begin{equation}
\label{hamil}
\begin{split}
&\dot{x}(t)=\frac{\partial H}{\partial
  \lambda}(t,x(t),\lambda(t),\lambda^0,u(t)),\\
&\dot{\lambda}(t)=-\frac{\partial H}{\partial x}(t,x(t),\lambda(t),\lambda^0,u(t)),\\
&\frac{\partial H}{\partial u}(t,x(t),\lambda(t),\lambda^0,u(t))=0,
\end{split}
\end{equation}
where
$$H(t,x, \lambda,\lambda^0,u):=\sum_{i=1}^mu_i\langle
\lambda,f_i(x)\rangle+\lambda^0 \left( \frac{1}{2}
u^TU(x)u+\alpha(x)^Tu+ g(t,x) \right)$$ 
is the {\em Hamiltonian} of the system.
An {\em extremal} is a $4$-tuple
$(x(\cdot),\lambda(\cdot),\lambda^0,u(\cdot))$ solution 
of the system of equations (\ref{hamil}). The extremal is said to be
{\em normal} if $\lambda^0\neq 0$ and {\em abnormal} if $\lambda^0=0$.

The relevance of singular trajectories in optimal control lies in the
fact that they are exactly the projections of abnormal extremals.
Note that a singular trajectory may be the projection of several
abnormal extremals, and also of a normal extremal.
A singular trajectory is said to be {\em strictly abnormal} if it is
not the projection of a normal extremal.
Notice that a singular trajectory is of corank one if and only if it
admits a unique (up to scalar normalization) abnormal extremal lift;
it is strictly abnormal and of corank one if and only if it admits a
unique extremal lift which is abnormal.

For a normal extremal, it is standard to adopt the normalization
$\lambda^0=-1$, and to derive the control $u$ as the feedback function
of $(x,\lambda)$
\begin{equation}\label{controlenormal}
u(t) =  \begin{pmatrix}u_1(t)\\ \vdots \\ u_m(t)\end{pmatrix}
=  U(x(t))^{-1} \begin{pmatrix} h_1(t)-\alpha_1(x(t))
  \\ \vdots \\ h_m(t)-\alpha_m(x(t))\end{pmatrix},
\end{equation}
for every $t\in[0,T]$, where
$ h_i(t) := \langle \lambda(t), f_i(x(t)) \rangle$, for $i=1,\ldots,m$.
In particular, normal extremals are smooth on $[0,T]$.

For abnormal extremals, the situation is much more involved, since
equations (\ref{hamil}) do not provide directly an expression for
abnormal controls. Abnormal extremals may be nonsmooth, and it is not
always possible to determine an explicit expression for
singular controls. Indeed, it follows from (\ref{hamil}) that
\begin{equation}\label{hiegalzero}
h_i(\cdot)\equiv 0 \ \textrm{on}\ [0,T],\ i=1,\ldots,m,
\end{equation}
along every abnormal extremal. At that point, in order to compute the
singular control, one usually differentiates iteratively
(\ref{hiegalzero}) with respect to $t$, until the control appears
explicitly (in an affine way). To recover the control, an
invertibility property is then required, which may not hold in
general.

In this paper, we prove that, in a generic context, such an
invertibility property is obtained with a minimal number of
differentiations (cf.\ Theorem \ref{thmA}). This is the concept of {\em
  minimal order}, defined in Definition \ref{defordre}.
Here, genericity means that the $(m+1)$-tuple $(f_0,\ldots,f_m)$
belongs to an open and dense subset of the set of vector fields
equipped with the Whitney topology. The corank one property is also
proved to hold generically. We obtain similar results in the driftless
case (cf.\ Theorem \ref{srthmA}).

In a preliminary step for deriving the above results, we establish a
theorem of independent interest, asserting that any trajectory of a
generic control-affine system satisfies $\dot{x}=0$ almost everywhere
on the set where the vector fields are linearly dependent
(cf.\ Theorem \ref{thmrang}).

\medskip

When considering optimal control problems, singular minimizing
trajectories may exist, and play a major role, since
they are not dependent on the specific minimization problem.
The issue of such minimizing trajectories was already well known in
the classical theory of calculus of variations (see for instance
\cite{bliss,Young}) and proved to be a major focus, during the forties, when
the whole domain eventually developed into optimal control theory (cf
\cite{Bonnard}).
The optimality status of singular trajectories was chiefly investigated in
\cite{BK1,trelatCOCV} in relation to control-affine systems with $m=1$, in
\cite{AS,LS,Montgomery,trelatCOCV} regarding driftless systems with
$m=2$ and in \cite{AS2,Sar} for general nonlinear control systems.

In this paper, we prove that, for generic systems with $m\geq 2$ (and
$m\geq 3$ in the driftless case), and for a fixed cost
$C_{U,\alpha,g}$, there does not exist minimizing singular
trajectories (cf.\ Corollaries \ref{cor1112} and \ref{srcor1112}).
We also prove that, given a fixed system $(\Sigma)$, singular
trajectories are strictly abnormal, generically with respect to the
cost (\ref{coutintro}) (cf.\ Propositions \ref{propstrictB11} and
\ref{srpropstrictB}).
We then show how the abovementioned results can be used to derive
regularity results for the value function and in the theory of
Hamilton-Jacobi equations, which in turn have applications for
stabilization and motion planning. 

\medskip

This paper is organized as follows.
Section \ref{sec2} is devoted to the statement of the main results,
firstly in the control-affine case, and secondly in the driftless
case. The consequences are detailed in Section \ref{sec3}, and
proofs are provided in Section \ref{sec4}.


\section{Statement of the main results}\label{sec2}
Let $M$ be a smooth, $n$-dimensional manifold.
Throughout the paper, $VF(M)$ denotes the set of smooth vector fields
on $M$, endowed with the $C^\infty$ Whitney topology.

\subsection{Trajectories of control-affine systems}
Let $T$ be a positive real number. Consider the control-affine system
\begin{equation}
\label{affine}
\dot{x}(t)=f_0(x(t))+\sum_{i=1}^mu_i(t)f_i(x(t)),
\end{equation}
where $(f_0,\ldots,f_m)$ is an $(m+1)$-tuple of smooth
vector fields on $M$, and the set of admissible controls
$u=(u_1,\ldots,u_m)$ is an open subset of $L^\infty([0,T],\Omega)$.

For every trajectory $x(\cdot):=x(\cdot,x_0,u)$ of (\ref{affine}),
define $I_{\mathrm{dep}}(x(\cdot))$ as the closed subset of $[0,T]$
\begin{equation}\label{Idep}
I_{\mathrm{dep}}(x(\cdot)) := \{ t\in[0,T]\ \vert\
\mathrm{rank}\{f_0(x(t)),\ldots,f_m(x(t))\}  < m+1 \} .
\end{equation}

Note that, on the open subset of $\R^n$ where
$\mathrm{rank}\{f_0,\ldots,f_m\}  = m+1$, there is a one-to-one
correspondence between trajectories and controls.
In contrast, on $I_{dep}(x(\cdot))$, there is no uniqueness of the
control associated to $x(\cdot)$; in particular, $x(\cdot)$ may be
associated to both singular and nonsingular controls.
This fact emphasizes the following result, which describes, in a generic
context, trajectories on the subset of $\R^n$ where
$\mathrm{rank}\{f_0,\ldots,f_m\}  < m+1$.

\begin{theorem}\label{thmrang}
Let $m<n$ be a nonnegative integer.
There exists an open and dense subset $O_{m+1}$ of
$VF(M)^{m+1}$ so that, if the $(m+1)$-tuple $(f_0,\ldots,f_m)$
belongs to $O_{m+1}$, then every trajectory $x(\cdot)$ of the associated
control-affine system
$
\dot{x}=f_0(x)+\sum_{i=1}^mu_if_i(x)
$
verifies
\begin{equation}\label{xpointIdep}
\dot{x}(t)=0,\ \textrm{for a.e.}\ t\in I_{\mathrm{dep}}(x(\cdot)).
\end{equation}
In addition, for every integer $N$, the set  $O_{m+1}$ can be chosen
so that its complement has codimension greater than $N$.
\end{theorem}

\begin{remark}\label{remidep}
At the light of the previous result, one can choose the admissible
control $u$ on $I_{dep}(x(\cdot))$ such that, for every $t\in
I_{dep}(x(\cdot))$, $u(t)$ consists of any $m$-tuple
$(\alpha_1,\ldots,\alpha_m)$ so that
$$
f_0(x(t))+\sum_{i=1}^m \alpha_if_i(x(t))=0.
$$
In particular, on any subinterval of $I_{dep}(x(\cdot))$, the
trajectory $x(\cdot)$ is constant, and the control can be
chosen constant as well.
\end{remark}

\begin{remark}\label{remtrivial}
A trajectory $x(\cdot)$ is said to be {\em trivial} if it reduces to a
point; otherwise it is said to be {\em nontrivial}.
It is clear that, if $I_{dep}(x(\cdot))\neq [0,T]$, then
$\dot{x}(t)\neq 0$ for $t\notin I_{dep}(x(\cdot))$ and $x(\cdot)$ is
nontrivial.

Let $x(\cdot)$ be a trajectory  of a control-affine system
associated to an $(m+1)$-tuple of $O_{m+1}$. As a consequence of
Theorem \ref{thmrang}, $x(\cdot)$ is trivial if and only if
$I_{dep}(x(\cdot))=[0,T]$.
\end{remark}


\subsection{Singular trajectories}
Recall that a singular trajectory $x(\cdot)$ is the projection of an abnormal
extremal $(x(\cdot),\lambda(\cdot))$. 
For $t\in[0,T]$ and $i,j\in\{0,\ldots,m\}$, we define
$$
h_i(t):=\langle \lambda(t),f_i(x(t)\rangle, \qquad
h_{ij}(t):=\langle \lambda(t),[f_i,f_j](x(t))\rangle .
$$
Along an abnormal extremal, we have for every $t\in[0,T]$,
\begin{equation}\label{s-affine}
h_0(t) = \textrm{constant},\ \ h_i(t)=0, \ \ i=1,\ldots,m.
\end{equation}
Differentiating (\ref{s-affine}), one gets, almost everywhere on
$[0,T]$,
\begin{equation}\label{sss-affine}
h_{i0}(t)+\sum_{j=1}^m h_{ij}(t)u_j(t)=0,\quad i\in\{0,\ldots,m\}.
\end{equation}

\begin{definition}
Along an abnormal extremal $(x(\cdot),\lambda(\cdot),u(\cdot))$ of the system
(\ref{affine}), the {\em Goh matrix} $G(t)$
at time $t\in [0,T]$ is
the $m\times m$ skew-symmetric matrix given by
\begin{equation}\label{goh-aff}
G(t):=\big(h_{ij}(t)\big)_{1\leq i,j\leq m}.
\end{equation}
\end{definition}

Since $G(t)$ is skew-symmetric, $\mathrm{rank}\ G(t)$ is
even, and Equation (\ref{sss-affine}) rewrites, almost everywhere on
$[0,T]$,
\begin{equation}\label{Au=b}
G(t) u(t) = b(t),
\end{equation}
with $b(t):=-(h_{i0}(t))_{1\leq i\leq m}$.

Note that, if $G(t)$ is invertible, then $u(t)$ is uniquely determined
by Equation (\ref{Au=b}). This only occurs for $m$ even.

If $m$ is odd, $G(t)$ is never invertible. However, a similar
construction is derived as follows. Define
\begin{equation}\label{Gbar}
\overline{G}(t):=\big(h_{ij}(t)\big)_{0\leq i,j\leq m}.
\end{equation}
Since $\overline{G}(t)$ is skew-symmetric, the determinant of
$\overline{G}(t)$ is the square of a polynomial $\overline{P}(t)$ in
the $h_{ij}(t)$ with degree $(m+1)/2$, called the {\em Pfaffian} (see
of $\overline{G}(t)$ \cite{artin}).
 From Equation (\ref{sss-affine}), $\overline{G}(t)$ is not
invertible, and thus, along the
extremal, $\overline{P}(t)=0$. After differentiation, one gets, almost
everywhere on $[0,T]$,
\begin{equation}\label{sss-affine2}
\{\overline{P},h_0\}(t)+\sum_{i=1}^m
u_j(t)\{\overline{P},h_j\}(t)=0.
\end{equation}
Define the $(m+1)\times m$ matrix $\widetilde{G}(t)$ as
$G(t)$ augmented with the row
$(\{\overline{P},h_j\}(t))_{1\leq j\leq m}$, and the $(m+1)$-dimensional
vector $\widetilde{b}(t)$ as $b(t)$ augmented with the coefficient
$-\{\overline{P},h_0\}(t)$.
Then, from Equations (\ref{Au=b}) and (\ref{sss-affine2}), there
holds, almost everywhere on $[0,T]$,
\begin{equation}\label{Atildeu=b}
\widetilde{G}(t)u(t)=\widetilde{b}(t).
\end{equation}
If $\widetilde{G}(t)$ is of rank $m$, then $u(t)$ is uniquely determined
by Equation (\ref{Atildeu=b}).

These facts, combined with Remark \ref{remidep}, motivate the following
definition.

\begin{definition}\label{defordre}
If $m$ is even (resp.\ odd),
a singular trajectory $x(\cdot)$ is said to be {\em of minimal order}
if:
\begin{itemize}
\item[(i)] $\dot{x}(t)=0$, for almost every $t\in I_{dep}(x(\cdot))$;
\item[(ii)] it admits an abnormal extremal lift such that, for almost every
$t\in [0,T]\setminus I_{dep}$, $\mathrm{rank}\
G(t)=m$ if $m$ is even, resp., $\mathrm{rank}\ \widetilde{G}(t)=m$ if
$m$ is odd.
\end{itemize}
\end{definition}

On the opposite, for arbitrary $m$,
a singular trajectory is said to be a {\em Goh trajectory} if it
admits an abnormal extremal
lift along which the Goh matrix is identically equal to $0$.

\begin{theorem}\label{thmA}
Let $m<n$ be a positive integer.
There exists an open and dense subset $O_{m+1}$ of 
$VF(M)^{m+1}$ so that, if the $(m+1)$-tuple
$(f_0,\ldots,f_m)$ belongs to $O_{m+1}$, then every nontrivial
singular trajectory of the associated control-affine system
$
\dot{x}(t)=f_0(x(t))+\sum_{i=1}^mu_i(t)f_i(x(t)),
$
is of minimal order and of corank one.
In addition, for every integer $N$, the set  $O_{m+1}$ can be chosen
so that its complement has codimension greater than $N$.
\end{theorem}

\begin{corollary}\label{corGoh2}
With the notations of Theorem~\ref{thmA} and if $m\geq 2$,
there exists an open and dense subset $O_{m+1}$ of $VF(M)^{m+1}$ 
so that every control-affine system defined with an
$(m+1)$-tuple of $O_{m+1}$ does not admit nontrivial Goh singular
trajectories.
\end{corollary}


\subsection{Minimizing singular trajectories}\label{subsec2.3}
We keep here the notations of the previous sections.
Consider the control-affine system
\begin{equation}\label{CAsyst}
\dot{x}(t)=f_0(x(t))+\sum_{i=1}^mu_i(t)f_i(x(t)),
\end{equation}
and the quadratic cost given by
\begin{equation}\label{coutUg}
C_{U,g}(T,u)=\frac{1}{2}\int_0^T \Big( u(t)^TU(x(t))u(t) + g(t,x(t)) \Big) dt,
\end{equation}
where $U\in\mathcal{S}^+_m(M)$ and $g\in C^\infty(\R\times M)$. Here,
$\mathcal{S}^+_m(M)$ denotes the set of smooth mappings $x\mapsto U(x)$ on
$M$, taking values in the set $\mathcal{S}^+_m$ of $m\times m$ real
positive definite matrices.

For $x_0\in M$ and $T>0$, define the optimal control problem
\begin{equation}\label{opt}
\inf \{ C_{U,g}(T,u)\ \vert\ E_{x_0,T}(u)=x \}.
\end{equation}
We next state two sets of genericity results, depending whether the
cost or the control system is fixed.

\subsubsection{Genericity w.r.t.\ the control system, with a fixed cost}
$\ $

\begin{proposition}\label{propstrictA}
Fix $U\in\mathcal{S}^+_m(M)$ and $g\in C^\infty(\R\times M)$.
There exists an open and dense subset $O_{m+1}$ of
$VF(M)^{m+1}$ such that
every nontrivial singular trajectory of
a control-affine system defined by an $(m+1)$-tuple
of $O_{m+1}$ is strictly abnormal for the optimal
control problem (\ref{opt}).
\end{proposition}

Corollary~\ref{corGoh2} together with Proposition~\ref{propstrictA}
yields the next corollary.

\begin{corollary}\label{cor1112}
Fix $U\in\mathcal{S}^+_m(M)$ and $g\in C^\infty(\R\times M)$.
Let $m\geq 2$ be an integer.
There exists an open and dense subset $O_{m+1}$ of $VF(M)^{m+1}$ so that 
the optimal control problem
(\ref{opt}) defined with an
$(m+1)$-tuple of $O_{m+1}$ does not admit nontrivial minimizing singular
trajectories.
\end{corollary}

\begin{remark}
In both previous results, the set  $O_{m+1}$ can be chosen
so that its complement has arbitrary codimension.
\end{remark}


\subsubsection{Genericity w.r.t.\ the cost, with a fixed control system}
We endow $\mathcal{S}^+_m(M)$ with the Whitney topology.
An $(m+1)-$ tuple $(f_0,\ldots,f_m)$ of $VF(M)^{m+1}$ is said to verify the
{\em Lie Algebra Rank Condition} if the Lie algebra generated by
$f_0,\ldots,f_m$ is of dimension $n$ at every point of $M$.

\begin{proposition}\label{propstrictB}
Fix $(f_0,\ldots,f_m)\in VF(M)^{m+1}$ so that the
Lie Algebra Rank Condition is satisfied and the zero control $u\equiv
0$ is not singular. Let $g\in C^\infty(\R\times M)$. Then, there
exists an open and dense subset $\mathcal{A}_m$ of
$\mathcal{S}^+_m(M)$ such that every nontrivial singular trajectory of
the control-affine system associated to the $(m+1)$-tuple
$(f_0,\ldots,f_m)$ is strictly abnormal for the optimal
control problem (\ref{opt}) defined with $U\in \mathcal{A}_m$
and $g$.
\end{proposition}

Assuming that the zero control $u\equiv 0$ is not singular is
a necessary hypothesis. Indeed, the fact that a control $u$ is 
singular is a property of the sole
$(m+1)$-tuple $(f_0,\ldots,f_m)$ and is independent of the cost.
On the other hand, every trajectory $x:=x(\cdot,x_0,0)$ associated 
to the zero control
is always the projection of the normal extremal $(x(\cdot),
0,-1,0)$ of any optimal
control problem (\ref{opt}). As a consequence, if the zero 
control is singular, such a trajectory $x(\cdot,x_0,0)$ cannot be
strictly abnormal.  

In order to handle the case of a singular zero control, it is
therefore necessary to consider more general quadratic costs such as
\begin{equation}\label{coutUg11}
C_{U,\alpha,g}(T,u)=\int_0^T\Big(\frac{1}{2}
u(t)^TU(x(t))u(t)+\alpha(x(t))^Tu(t) 
+ g(t,x(t)) \Big) dt,
\end{equation}
where $U\in\mathcal{S}^+_m(M)$, $\alpha\in C^{\infty}(M,\R^m)$ 
and $g\in C^\infty(\R\times M)$.

\begin{proposition}\label{propstrictB11}
Fix $(f_0,\ldots,f_m)\in VF(M)^{m+1}$ satisfying the
Lie Algebra Rank Condition and $g\in C^\infty(\R\times M)$.
Then, there exists an open and dense subset $\mathcal{B}_m$ of
$\mathcal{S}^+_m(M)\times C^{\infty}(M,\R^m)$ such that every 
nontrivial singular trajectory of
the control-affine system associated to the $(m+1)$-tuple
$(f_0,\ldots,f_m)$ is strictly abnormal for the optimal
control problem (\ref{CAsyst})-(\ref{coutUg11}) defined with 
$(U,\alpha)\in \mathcal{B}_m$ and $g$.
\end{proposition}


\subsection{Driftless control-affine systems}

Let $T$ be a positive real number. Consider the driftless
control-affine system 
\begin{equation}
\label{sraffine}
\dot{x}(t)=\sum_{i=1}^mu_i(t)f_i(x(t)),
\end{equation}
where $(f_1,\ldots,f_m)$ is an $m$-tuple of smooth
vector fields on $M$, and the set of admissible controls
$u=(u_1,\ldots,u_m)$ is an open subset of $L^\infty([0,T],\Omega)$.

For every trajectory $x(\cdot):=x(\cdot,x_0,u)$ of (\ref{affine}),
define $I_{\mathrm{dep}}(x(\cdot))$ as the closed subset of $[0,T]$
$$
I_{\mathrm{dep}}(x(\cdot)) := \{ t\in[0,T]\ \vert\
\mathrm{rank}\{f_1(x(t)),\ldots,f_m(x(t))\}  < m \} .
$$

\begin{theorem}\label{srthmrang}
Let $m\leq n$ be a positive integer.
There exists an open and dense subset $O_{m}$ of
$VF(M)^{m}$ so that, if the $m$-tuple $(f_1,\ldots,f_m)$
belongs to $O_{m}$, then every trajectory $x(\cdot)$ of the associated
driftless control-affine system
$
\dot{x}=\sum_{i=1}^m u_i f_i(x)
$
verifies 
$$
\dot{x}(t)=0,\ \textrm{for a.e.}\ t\in I_{\mathrm{dep}}(x(\cdot)).
$$
In addition, for every integer $N$, the set  $O_{m}$ can be chosen
so that its complement has codimension greater than $N$.
\end{theorem}

\subsubsection{Singular trajectories}
Let $x(\cdot)$ be a singular trajectory, projection of an abnormal
extremal $(x(\cdot),\lambda(\cdot))$. Similarly to the previous section,
we define, for $t\in[0,T]$ and $i,j\in\{1,\ldots,m\}$,
$$
h_i(t):=\langle \lambda(t),f_i(x(t)\rangle, \qquad
h_{ij}(t):=\langle \lambda(t),[f_i,f_j](x(t))\rangle .
$$
For every $t\in[0,T]$,
\begin{equation}\label{srs-affine}
h_i(t)=0, \ \ i=1,\ldots,m.
\end{equation}
Differentiating (\ref{srs-affine}), one gets, almost everywhere on
$[0,T]$,
\begin{equation}\label{srsss-affine}
\sum_{j=1}^m h_{ij}(t)u_j(t)=0,\quad i\in\{1,\ldots,m\}.
\end{equation}

\begin{definition}
Along an abnormal extremal $(x(\cdot),\lambda(\cdot),u(\cdot))$ of the system
(\ref{affine}), the {\em Goh matrix} $G(t)$
at time $t\in [0,T]$ is
the $m\times m$ skew-symmetric matrix given by
\begin{equation}\label{srgoh-aff}
G(t):=\big(h_{ij}(t)\big)_{1\leq i,j\leq m}.
\end{equation}
\end{definition}

Since $G(t)$ is skew-symmetric, $\mathrm{rank}\ G(t)$ is
even, and Equation (\ref{srsss-affine}) rewrites, almost everywhere on
$[0,T]$,
\begin{equation}\label{srAu=b}
G(t) u(t) = 0.
\end{equation}

Note that, if $\mathrm{rank}\ G(t)=m-1$, 
one can deduce from~(\ref{srAu=b})  an
expression for $u(t)$, up to time reparameterization.
This only occurs for $m$ odd.

If $m$ is even, $\mathrm{rank}\ G(t)$ is always smaller than
$m-1$. However, a similar 
construction is derived as follows. The determinant of
$G(t)$ is the square of the  Pfaffian $P(t)$, and, along the
extremal, $P(t) \equiv 0$. After differentiation, one gets, almost
everywhere on $[0,T]$,
\begin{equation}\label{srsss-affine2}
\sum_{i=1}^m
u_j(t)\{ P,h_j\}(t)=0.
\end{equation}
Define the $(m+1)\times m$ matrix $\widetilde{G}(t)$ as
$G(t)$ augmented with the row
$(\{ P,h_j\}(t))_{1\leq j\leq m}$.
Then, from Equations (\ref{srAu=b}) and (\ref{srsss-affine2}), there
holds, almost everywhere on $[0,T]$,
\begin{equation}\label{srAtildeu=b}
\widetilde{G}(t)u(t)=0.
\end{equation}
If $\widetilde{G}(t)$ is of rank $m-1$, one can deduce
from~(\ref{srAtildeu=b})  an 
expression for $u(t)$, up to time reparameterization.

\begin{definition}\label{srdefordre}
If $m$ is odd (resp.\ even), a singular trajectory $x(\cdot)$ is said
to be {\em of minimal order} 
if:
\begin{itemize}
\item[(i)] $\dot{x}(t)=0$, for almost every $t\in I_{dep}(x(\cdot))$;
\item[(ii)] it admits an abnormal extremal lift such that, for almost every
$t\in [0,T]\setminus I_{dep}$, $\mathrm{rank}\
G(t)=m-1$ if $m$ is odd, resp., $\mathrm{rank}\ \widetilde{G}(t)=m-1$ if
$m$ is even.
\end{itemize}
\end{definition}

On the opposite, for arbitrary $m$,
a singular trajectory is said to be a {\em Goh trajectory} if it
admits an abnormal extremal
lift along which the Goh matrix is identically equal to $0$.

\begin{theorem}\label{srthmA}
Let $m$ be an integer such that $2 \leq m \leq n$.
There exists an open and dense subset $O_{m}$ of $VF(M)^{m}$ so that, 
if the $m$-tuple $(f_1,\ldots,f_m)$ belongs to $O_{m}$, then every nontrivial
singular trajectory of the associated driftless control-affine system
$\dot{x}(t)=\sum_{i=1}^mu_i(t)f_i(x(t))$
is of minimal order and of corank one.
In addition, for every integer $N$, the set  $O_{m}$ can be chosen
so that its complement has codimension greater than $N$.
\end{theorem}

\begin{corollary}\label{srcorGoh2}
With the notations of Theorem~\ref{srthmA} and if $m\geq 3$,
there exists an open and dense subset $O_{m}$ of $VF(M)^{m}$ so that 
every driftless control-affine system defined with an
$m$-tuple of $O_{m}$ does not admit nontrivial Goh singular
trajectories.
\end{corollary}


\subsubsection{Minimizing singular trajectories}\label{subsubsec2.4.2}
Consider the optimal control problem
associated to the driftless control-affine system
\begin{equation}\label{srCAsyst}
\dot{x}(t)=\sum_{i=1}^mu_i(t)f_i(x(t)),
\end{equation}
with the quadratic cost given by
\begin{equation}\label{srcoutUg}
C_{U,g}(T,u)=\frac{1}{2}\int_0^T \Big( u(t)^TU(x(t))u(t) + g(t,x(t))
\Big) dt, 
\end{equation}
where $U\in\mathcal{S}^+_m(M)$ and $g\in C^\infty(\R\times M)$.

For $x_0\in M$ and $T>0$, define the optimal control problem
\begin{equation}\label{sropt}
\inf \{ C_{U,g}(T,u)\ \vert\ E_{x_0,T}(u)=x \}.
\end{equation}
We next state genericity results with respect to the control system, with a
fixed cost.

\begin{proposition}\label{srpropstrictA}
Fix $U\in\mathcal{S}^+_m(M)$ and $g\in C^\infty(\R\times M)$.
There exists an open and dense subset $O_{m}$ of
$VF(M)^{m}$ such that every nontrivial singular trajectory of
a driftless control-affine system defined by an $m$-tuple
of $O_{m}$ is strictly abnormal for the optimal
control problem (\ref{sropt}).
\end{proposition}

Corollary~\ref{srcorGoh2} together with Proposition~\ref{srpropstrictA}
yields the next corollary.

\begin{corollary}\label{srcor1112}
Fix $U\in\mathcal{S}^+_m(M)$ and $g\in C^\infty(\R\times M)$.
Let $m\geq 3$ be an integer.
There exists an open and densee subset $O_{m}$ of $VF(M)^{m}$ so that the optimal control problem
(\ref{sropt}) defined with an
$m$-tuple of $O_{m}$ does not admit nontrivial minimizing singular
trajectories.
\end{corollary}

\begin{remark}
In both previous results, the set  $O_{m}$ can be chosen
so that its complement has arbitrary codimension.
\end{remark}

We also have have a genericity result with respect to the cost, with a
fixed control system.

\begin{proposition}\label{srpropstrictB}
Fix $(f_1,\ldots,f_m)\in VF(M)^{m}$ so that the
Lie Algebra Rank Condition is satisfied. Let $g\in C^\infty(\R\times
M)$. Then, there 
exists an open and dense subset $\mathcal{A}_m$ of
$\mathcal{S}^+_m(M)$ such that every nontrivial singular trajectory of
the driftless control-affine system associated to the $m$-tuple
$(f_1,\ldots,f_m)$ is strictly abnormal for the optimal
control problem (\ref{sropt}) defined with $U\in \mathcal{A}_m$
and $g$.
\end{proposition}

\begin{remark}
In the driftless case, the control $u\equiv 0$ is always singular, but
corresponds to a trivial trajectory. Therefore, in opposition to the
control-affine case, it is not necessary to add the linear term
$\alpha(x)^Tu$ in the cost.
\end{remark}


\section{Consequences}  
\label{sec3}
\subsection{Regularity of the value function}
Consider the optimal control problem (\ref{opt}), associated to the
control-affine system (\ref{CAsyst}) and the cost (\ref{coutUg}). The value
function is defined by
\begin{equation}\label{opt1}
S_{x_0,T}(x) := \inf \{ C_{U,g}(T,u)\ \vert\ E_{x_0,T}(u)=x \},
\end{equation}
for every $x\in\R^n$ (with, as usual,
$\inf\emptyset:=-\infty$). We assume in the sequel that all data are
analytic.

The regularity of $S_{x_0,T}$ is closely related to the existence of
nontrivial minimizing singular trajectories starting from $x_0$.
It is proved in \cite{trelatJDCS} that, in the absence of minimizing
singular trajectories, the value function is continuous and
subanalytic (see e.g.\ \cite{Hironaka} for a definition of a subanalytic
function).
For driftless control-affine systems and $g\equiv 0$, the value
function coincides with the square of a sub-Riemannian distance
(see \cite{Be} for an introduction to sub-Riemannian geometry). In
particular, in this case, the value function is always continuous, but
the trivial trajectory $x(\cdot)\equiv x_0$ is always minimizing and
singular. Moreover, if there is no nontrivial minimizing singular
trajectories, then the value function is subanalytic outside $x_0$
(see \cite{Ag,AgS}). This situation holds for generic distributions of
rank greater than or equal to three (see \cite{AG,CJT}).

The results of Section \ref{subsec2.3} have the following
consequence on the regularity of $S_{x_0,T}$.

\begin{corollary}
With the notations of Corollary \ref{cor1112}, and
if in addition the functions $g$, $U$, and
the vector fields of the $(m+1)$-tuple in $O_{m+1}$
are analytic, then the associated value function $S_{x_0,T}$ is continuous
and subanalytic on its domain of definition.
\end{corollary}

\begin{remark}
If there exists a nontrivial minimizing trajectory, the value function
may fail to be subanalytic or even continuous.
For example, consider the control-affine system
in $\mathbb{R}^2$ given by
\begin{equation*}
\dot{x}(t)=1+y(t)^2,\quad
\dot{y}(t)=u(t),
\end{equation*}
and the cost
$C(T,u)=\int_0^T u(t)^2dt$.
The trajectory $(x(t)=t,y(t)=0)$, associated
to the control $u=0$, is a nontrivial minimizing singular trajectory,
and the value function $S_{(0,0),T}$ has the asymptotic expansion,
near the point $(T,0)$,
\begin{equation*}
S_{(0,0),T}(x,y)= \frac{1}{4}\frac{y^4}{x-T} + 
\frac{y^4}{x-T}\,\mathrm{exp}\left(-\frac{y^2}{x-T}\right) 
+\mathrm{o}\left(\frac{y^4}{x-T}\,\mathrm{exp}\left(-\frac{y^2}{x-T}\right)\right)
\end{equation*}
(see \cite{trelatJDCS} for details).
Hence, it is not continuous, nor subanalytic, at the point $(T,0)$.
\end{remark}

In the driftless control-affine case, by using results of Section
\ref{subsubsec2.4.2}, we derive the following similar consequence.

\begin{corollary}
With the notations of Corollary \ref{srcor1112}, and
if in addition the functions $g$, $U$, and
the vector fields of the $m$-tuple in $O_{m}$
are analytic, then the associated value function $S_{x_0,T}$ is
subanalytic outside $x_0$.
\end{corollary}


\subsection{Regularity of viscosity solutions of Hamilton-Jacobi equations}
Assume that the assumptions of the previous subsection hold.
It is standard (see \cite{Evans,Lions}) that the value function
$v(t,x)=S_{x_0,t}(x)$ is a viscosity solution of the Hamilton-Jacobi
equation
\begin{equation}\label{hjbeq}
\frac{\partial v}{\partial t} + \mathcal{H}\left(x,\frac{\partial v}{\partial
    x}\right)= g(t,x),
\end{equation}
where
$ 
\mathcal{H}(x,p) = 
\langle p,f_0(x)\rangle + \frac{1}{2}\sum_{i,j=1}^m
(U^{-1}(x))_{ij} \langle p,f_i(x)\rangle \langle p,f_j(x)\rangle  $.

Conversely, the viscosity solution of (\ref{hjbeq}) with analytic
Dirichlet-type conditions is subanalytic, as soon as the corresponding
optimal control does not admit minimizing singular trajectories
(see \cite{trelatHJB}). Using the results of the previous sections,
this situation holds generically if $m\geq 2$ (and, similarly for
driftless control-affine systems, if $m\geq 3$).

As a consequence, the analytic singular set $\mathrm{Sing}(v)$ of the
viscosity solution $v$, i.e., the subset of $\R^n$ where $v$ is not
analytic, is a (subanalytic) stratified manifold of codimension
greater than or equal to one (see \cite{Tamm} for more details on the
subject). Since $\mathrm{Sing}(v)$ is also the locus
where characteristic curves intersect, the abovementioned property
turns out to be instrumental for the global convergence of numerical
schemes for Equation (\ref{hjbeq}) (see \cite{Evans}). Indeed, the
analytic singular set must be as ``nice'' as possible in order to
integrate energy functions on the set of characteristic curves.


\subsection{Applications to stabilization and motion planning}
For a driftless control-affine system verifying the Lie
Algebra Rank Condition, there exist general stabilizing strategies
stemming from dynamic programming. As usual, the stabilizing feedback
is computed using the gradient of the value function $S$ for a suitable
optimal control problem. Of course this is only possible outside the
singular set $\mathrm{Sing}(S)$, and one must device another
construction for the feedback on $\mathrm{Sing}(S)$.
Let us mention two such strategies, the first one providing an hybrid
feedback (see \cite{PrieurTrelatSIAM}), and the second one a
smooth repulsive stabilizing (SRS) feedback (see \cite{Rifford1,Rifford2}).
Both strategies crucially rely on the fact that $\mathrm{Sing}(S)$ is
a stratified manifold of codimension greater than or equal to one.

As seen before, the latter fact holds generically in the analytic
category for $m\geq 3$.

\medskip

On the other hand, the absence of singular minimizing trajectories is
the basic requirement for the convergence of usual algorithms in
optimal control (such as direct or indirect methods, see
e.g. \cite{Betts,Pesch}). We have proved that
this situation holds generically for control-affine systems if $m\geq
2$, and for driftless control-affine systems if $m\geq 3$.

As a final application,
consider a driftless control-affine system verifying the Lie
Algebra Rank Condition. According to Proposition \ref{srpropstrictB}, it is
possible to choose a (generic) cost function $C_{U,g}$ such that all
singular trajectories are strictly abnormal. Combining that fact with
\cite[Theorem 1.1]{RiffordTrelatMA}, we deduce that there exists a
dense subset $N$ of $\R^n$ such that every point of $N$ is reached by
a unique minimizing trajectory, which is moreover nonsingular.
As a consequence, a shooting method with target in $N$ will converge.
That fact may be used for solving (at least approximately) motion
planning problems.


\section{Proofs of the results}
\label{sec4}
\subsection{Proof of Theorem~\ref{thmrang} and of Theorem~\ref{srthmrang}}
Every trajectory of the control-affine system
$
\dot{x}=f_0(x)+\sum_{i=1}^mu_if_i(x)
$
is also a trajectory of the driftless control system
$
\dot{x}=\sum_{i=0}^mu_if_i(x),
$
with $u_0\equiv 1$. Therefore, Theorem~\ref{thmrang} follows from
Theorem~\ref{srthmrang}, whose proof is provided next. 

Let $x(\cdot)=x(\cdot,x_0,u)$ be a trajectory of the
driftless control system
$\dot{x}=\sum_{i=1}^m u_if_i(x)$, with $2\leq m\leq n$.
Consider the set $I_{\mathrm{dep}}(x(\cdot))$ defined by
(\ref{Idep}). We argue by contraposition, and assume that
$I_{\mathrm{dep}}(x(\cdot))$ contains a subset $I$ of positive measure
such that $\dot{x}(t)\neq 0$ for $t\in I$.
Since Lebesgue points of $u$ are of full Lebesgue measure, we assume
that $u$ is continuous on $I$.

Up to considering a subset of $I$, and relabeling the $f_i$'s, we
assume that, for every $t\in I$:
\begin{itemize}
\item[(i)] there exists $1\leq k<m$ such that
$$\mathrm{rank}\{f_1(x(t)),\ldots,f_m(x(t))\}=k;$$
\item[(ii)] $f_1(x(t)),\ldots,f_k(x(t))$ are linearly independent, and
  thus, there exist real numbers $\alpha_i^j(t)$, $i=1,\ldots,k$,
  $j=k+1,\ldots,m$, such that
\begin{equation*}
f_j(x(t)) = \sum_{i=1}^k \alpha_i^j(t) f_i(x(t)),\quad j=k+1,\ldots,m.
\end{equation*}
Therefore,
$\displaystyle
\dot{x}(t) = \sum_{i=1}^k \delta_i(t) f_i(x(t)),
$
where $\displaystyle\delta_i(t):=u_i(t)+\sum_{j=k+1}^m \alpha_i^j(t)u_j(t)$;
\item[(iii)] $\delta_1(t)\neq 0$.
\end{itemize}

\begin{remark}\label{remrangouvert}
Up to reducing $I$, we furthermore assume that $I$ is contained in an
open interval $\mathcal{I}$ on which
$\mathrm{rank}\{f_1(x(t)),\ldots,f_k(x(t))\}=k$.
\end{remark}

Set $\mathrm{ad}^0g(h)=h$, where $g,h\in\mathrm{VF}(M)$,
and $\mathrm{ad}^kg(h)=[g,\mathrm{ad}^{k-1}g(h)]$ for $k\geq 1$.
The length of the iterated Lie bracket $[f_{i_1},[f_{i_2},[\cdots,
f_{i_k}]\cdots]$ of $f_1,\ldots,f_m$ is the integer $k$.

\begin{proposition}\label{lemroutine}
Let $N$ be a positive integer. There exists a subset $J_N\subset I$ of
positive measure such that, for every $t\in J_N$, and every
$\ell\in\{1,\ldots,N\}$,
\begin{equation}\label{rout0}
\delta_1(t)^{\ell-1} \mathrm{ad}^{\ell-1}f_1(f_m) (x(t)) =
h_t^\ell(x(t)) + R_t^\ell(x(t)),
\end{equation}
where
\begin{itemize}
\item $h_t^\ell(x(t))\in\mathrm{Span}\{f_1(x(t)),\ldots,f_k(x(t))\}$,
\item $R_t^\ell$ is a linear combination of iterated Lie
  brackets of $f_1,\ldots,f_m$, of length smaller than $\ell-1$, and
  of iterated Lie brackets of $f_1,\ldots,f_k$, of length smaller than
  or equal to $\ell$.
\end{itemize}
\end{proposition}

\begin{proof}
For $t\in I$, let $F_t\in VF(M)$ be the vector field defined by
$$
F_t(x):=\sum_{i=1}^k \delta_i(t) f_i(x).
$$
Notice that $\dot{x}(t)=F_t(x(t))$, for $t\in I$.
For the argument of Proposition~\ref{lemroutine}, we need the
following lemma.

\begin{lemma}\label{lem-recursif}
Consider a set $J\subset I$ of positive measure and $h\in VF(M)$ so
that
$h(x(t))\in \mathrm{Span}\{f_1(x(t)),\ldots,f_k(x(t))\}$ on
$J$,
i.e., for every $t\in J$, there exist real numbers $\beta_i(t)$, 
$i=1,\ldots,k$, such that
\begin{equation}\label{h00}
h(x(t)) = \sum_{i=1}^k \beta_i(t) f_i(x(t)).
\end{equation}
For $t\in J$, define $g_t\in VF(M)$ by 
$$
g_t(x):= h(x)- \sum_{i=1}^k \beta_i(t) f_i(x).
$$
Then, there exists a set $J'\subset J$ of positive measure such
that
\begin{equation}\label{eq00}
[F_t,g_t](x(t))\in \mathrm{Span}\{f_1(x(t)),\ldots,f_k(x(t))\} \quad \hbox{ on
}J'.
\end{equation}
\end{lemma}

\begin{proof}[Proof of Lemma~\ref{lem-recursif}]
 Using Remark \ref{remrangouvert},
there exist $e_j\in VF(M)$, $k+1\leq j\leq n$, so
that, for every $t\in \mathcal{I}$, the vectors
$f_1(x(t))$, $\dots$, $f_k(x(t))$, $e_{k+1}(x(t))$, $\dots$, $e_n(x(t))$ span
$T_{x(t)}M$. Thus,
there exist $n$ smooth functions $b_i$, $1\leq
i\leq n$, defined on $M$, such that
$$
h(x)=\sum_{i=1}^kb_i(x)f_i(x)+\sum_{i=k+1}^nb_i(x)e_i(x),
$$
for $x$ in an open neighborhood of $x(\mathcal{I})$.
For $i=1,\dots, n$, define $\beta_i(t):=b_i(x(t))$, for $t\in \mathcal{I}$
(this notation is consistent with (\ref{h00})).
The $\beta_i's$ are absolutely continuous on $\mathcal{I}$ and
differentiable everywhere on $J$. For $i=k+1,\dots, n$, there holds
$\beta_i\equiv 0$ on $J$ and therefore, it follows that
$\dot\beta_i\equiv 0$ on a subset $J'\subset J$ of full
measure (cf. \cite[Lemma p.~177]{Ru}).

For $t\in J$, using that $g_t(x(t))=0$, and $F_t(x(t))=\dot{x}(t)$, it
holds
\begin{eqnarray*}
[F_t,g_t](x(t)) &=& dg_t\circ F_t(x(t)) \\
&=& \sum_{i=1}^k (db_i(x(t)).\dot{x}(t)) f_i(x(t)) +
 \sum_{i=k+1}^n (db_i(x(t)).\dot{x}(t)) e_i(x(t)) \\
&=& \sum_{i=1}^k \dot{\beta}_i(t) f_i(x(t)) +
 \sum_{i=k+1}^n \dot{\beta}_i(t) e_i(x(t)).
\end{eqnarray*}
On $J'$, the second sum of the right-hand side of the last equation
vanishes, and the lemma follows.
\end{proof}

Applying Lemma \ref{lem-recursif} to $h=f_m$ and $J=I$, we get
\begin{equation*}
[F_t,g_t^1](x(t))\in \mathrm{Span}\{f_1(x(t)),\ldots,f_k(x(t))\} \quad \hbox{ on
}J_1,
\end{equation*}
where $J_1\subset I$ and
$g_t^1 := f_m-\sum_{i=1}^k \alpha_i^m(t)f_i$.

Set $h_t^1=[F_t,g_t^1]$. We next iterate the above procedure, for
$1\leq\ell\leq N$. Assume that the vector fields $h_t^\ell$,
$g_t^\ell$, and the set $J_\ell$ of positive measure are defined, such
that
$h_t^\ell(x(t))\in \mathrm{Span}\{f_1(x(t)),\ldots,f_k(x(t))\}$ on
$J_\ell$.
For every $t\in J_\ell$, let $\beta_i^\ell(t)$,
$i=1,\ldots,k$, be the real numbers such that
\begin{equation*}
h_t^\ell(x(t)) = \sum_{i=1}^k \beta_i^\ell(t) f_i(x(t)),
\end{equation*}
and define $g_t^{\ell+1}\in VF(M)$ by 
$
g_t^{\ell+1}:= h_t^\ell- \sum_{i=1}^k \beta_i^\ell(t) f_i.
$
Set $ h_t^{\ell+1} := [F_t,g_t^{\ell+1}]$.
Applying Lemma \ref{lem-recursif}, there exists a subset
$J_{\ell+1}\subset J_\ell$ of positive measure such that
$h_t^{\ell+1}(x(t))\in \mathrm{Span}\{f_1(x(t)),\ldots,f_k(x(t))\}$
on $J_{\ell+1}$.

For $t\in J_N$, and for $\ell=1,\ldots,N$, we express $h_t^\ell(x(t))$
using iterated Lie brackets of $f_1,\ldots,f_m$, and an easy induction
yields (\ref{rout0}).
\end{proof}

Combining Proposition \ref{lemroutine} with routine transversality
arguments (see for instance \cite{BK} and \cite{CJT}), it follows that
the $(m+1)$-tuple $(f_0,\ldots,f_m)$ belongs to a closed subset of
$VF(M)^{m+1}$ of codimension greater than or equal to
$N$. Theorem~\ref{srthmrang} follows. 

\begin{remark}\label{remtransverse}
The fact that $f_1(x(t))\neq 0$ is essential in order
to derive, from (\ref{rout0}), an infinite number of independent
relations, and then to apply the abovementioned transversality
arguments.
\end{remark}


\subsection{Proof of Theorem \ref{thmA}}
The minimal order and corank one properties are proved separately in
the following lemmas.

\begin{lemma}\label{keylemma1}
There exists an open and dense subset $O^1_{m+1}$ of $VF(M)^{m+1}$ so that, 
if the $(m+1)$-tuple $(f_0,\ldots,f_m)$ belongs to $O^1_{m+1}$, then every
singular trajectory of the associated control-affine system
$
\dot{x}(t)=f_0(x(t))+\sum_{i=1}^mu_i(t)f_i(x(t))
$
is of minimal order.
In addition, for every integer $N$, the set $O^1_{m+1}$ can be chosen
so that its complement has codimension greater than $N$.
\end{lemma}

\begin{lemma}\label{keylemma2}
There exists an open and dense subset $O_{m+1}$ of $O^1_{m+1}$ so that, 
if the $(m+1)$-tuple $(f_0,\ldots,f_m)$ belongs to $O_{m+1}$, then every nontrivial
singular trajectory of the associated control-affine system
$
\dot{x}(t)=f_0(x(t))+\sum_{i=1}^mu_i(t)f_i(x(t))
$
is of corank one.
In addition, for every integer $N$, the set  $O^2_{m+1}$ can be chosen
so that its complement has codimension greater than $N$.
\end{lemma}

The conclusion of Theorem \ref{thmA} with the set $O_{m+1}$ whose
existence is stated above.

\subsubsection{Proof of Lemma \ref{keylemma1}}
 From Theorem \ref{thmrang}, there exists an open and dense subset $O_{m+1}^{11}$
of $VF(M)^{m+1}$ such that, if $(f_0,\ldots,f_m)\in O_{m+1}^{11}$, then
every trajectory $x(\cdot)$ of $\dot{x}=f_0(x)+\sum_{i=1}^m u_if_i(x)$
verifies Item (i) of Definition \ref{defordre}.

It is therefore enough to show the existence of an open and dense subset
$O_{m+1}^{12}$ of $VF(M)^{m+1}$ such that, if $(f_0,\ldots,f_m)\in
O_{m+1}^{12}$, then every singular trajectory $x(\cdot)$ of
$\dot{x}=f_0(x)+\sum_{i=1}^m u_if_i(x)$ verifies Item (ii) of
Definition \ref{defordre}.
Then, by choosing $O^1_{m+1}:=O_{m+1}^{11}\cap O_{m+1}^{12}$, the
conclusion of 
Lemma~\ref{keylemma1} follows.

Consider a singular trajectory $x(\cdot):=x(\cdot,x_0,u)$ of
$\dot{x}=f_0(x)+\sum_{i=1}^m u_if_i(x)$, admitting
an abnormal extremal $(x(\cdot),\lambda(\cdot))$. Assume that there
exists $J\subset [0,T]\setminus I_{dep}(x(\cdot))$ of
positive measure such that
$G(t)$ is not of rank $m$ if $m$ is even, respectively,
$\widetilde{G}(t)$ is not of rank $m$ if $m$ is odd.
We will show that the $(m+1)$-tuple  $(f_0,\ldots,f_m)$ belongs to a
subset of arbitrary codimension in $VF(M)^{m+1}$ whose complement
contains an open and dense subset.

Note that, on $[0,T]\setminus I_{dep}(x(\cdot))$, the vector fields
$f_0(x(t)),\ldots, f_m(x(t))$ are linearly independent.
The remaining part of the argument consists of reformulating the
problem in order to follow the chain of arguments in the proof of
\cite[Theorem 2.4]{CJT}
concerning the case of everywhere linearly independent vector fields.
For that purpose, we distinguish the cases $m$ even and $m$ odd.

Assume first that $m$ is even. As in (\ref{Gbar}), define, for $t\in
J$,
$
\overline{G}(t):=\big(h_{ij}(t)\big)_{0\leq i,j\leq m}.
$
 From (\ref{Au=b}), we have, for $t\in J$,
$$ 
\overline{G}(t) = 
\begin{pmatrix} 
  0 & \big(G(t)u(t)\big)^T \\ 
  -G(t)u(t) & G(t) 
\end{pmatrix} .
$$
Since the ranks of both $\overline{G}(t)$ and $G(t)$ are even, 
they must be equal, for $t\in J$, and hence,
the rank of $\overline{G}(t)$ is smaller than $m$ on $J$.
This is exactly the starting point of the proof of
\cite[Lemma 3.8]{CJT}.
The machinery of \cite{CJT} then applies and we deduce that the
$(m+1)$-tuple  $(f_0,\ldots,f_m)$ belongs to a subset of arbitrary
codimension in $VF(M)^{m+1}$ whose complement contains an open and
dense subset $O_{m+1}^2$ of $VF(M)^{m+1}$.

Assume next that $m$ is odd.
Define the $(m+2)\times (m+1)$ matrix $\widehat{G}(t)$ as
$\overline{G}(t)$ augmented in the last row with
$(\{\overline{P},h_j\}(t))_{0\leq j\leq m}$.

\begin{lemma}\label{lemrank}
With the notations above,
$
\mathrm{rank}\ \widehat{G}(t) \leq \mathrm{rank}\ \widetilde{G}(t)+1.
$
\end{lemma}

\begin{proof}
It amounts to show that, $\xi\in\ker\widetilde{G}(t)$ implies 
$(0,\xi)\in\ker\widehat{G}(t)$. This follows from the fact that
if $\widetilde{G}(t)\xi=0$, then $G(t)\xi=0$, and thus $\xi$ is
orthogonal to the range of $G(t)$ since $G(t)$ is skew-symmetric.
\end{proof}

Using Lemma \ref{lemrank}, the rank of $\widehat{G}(t)$ is less than
$m+1$ on $J$.
This is exactly the starting point of the proof of \cite[Lemma
3.9]{CJT}. The machinery of \cite{CJT} then applies and we deduce that
the $(m+1)$-tuple  $(f_0,\ldots,f_m)$ belongs to a subset of arbitrary
codimension in $VF(M)^{m+1}$ whose complement contains an open and
dense subset $O_{m+1}^{12}$ of $VF(M)^{m+1}$.


\subsubsection{Proof of Lemma \ref{keylemma2}}
We argue by contraposition.
Consider a nontrivial singular trajectory $x(\cdot):=x(\cdot,x_0,u)$ of
$\dot{x}=f_0(x)+\sum_{i=1}^m u_if_i(x)$, with $(f_0,\ldots,f_m)\in
O_{m+1}^1$. Assume that $x(\cdot)$ admits two abnormal extremal
lifts $(x(\cdot),\lambda^{[1]}(\cdot))$ and
$(x(\cdot),\lambda^{[2]}(\cdot))$ such that, for some $t_0 \in [0,T]$,
$\lambda^{[1]}(t_0)$ and $\lambda^{[2]}(t_0)$ are linearly independent.
By linearity, $\lambda^{[1]}(\cdot)$ and $\lambda^{[2]}(\cdot)$ are linearly independent
everywhere on $[0,T]$.
Since $x(\cdot)$ is nontrivial, it follows from
Remark~\ref{remtrivial} that there exists a nonempty subinterval $J$
of $[0,T]\setminus I_{dep}(x(\cdot))$.
We are now in a position to exactly follow the arguments of \cite{CJT}
corresponding to the corank one property, i.e., \cite[Lemma 4.4]{CJT}.


\subsection{Proof of Theorem~\ref{srthmA}}
We start with the proof of the statement dealing with the minimal
order property. 

 From Theorem~\ref{srthmrang}, there exists an open and dense subset 
$O_{m}^{1}$ of $VF(M)^{m}$ such that, if $(f_1,\ldots,f_m)\in
O_{m}^{1}$, then
every trajectory $x(\cdot)$ of $\dot{x}=\sum_{i=1}^m u_if_i(x)$
verifies Item (i) of Definition~\ref{srdefordre}.

It is therefore enough to show the existence of an open and dense subset
$O_{m}^{2}$ of $VF(M)^{m}$ such that, if $(f_1,\ldots,f_m)\in
O_{m}^{2}$, then every singular trajectory $x(\cdot)$ of
$\dot{x}=\sum_{i=1}^m u_if_i(x)$ verifies Item (ii) of
Definition~\ref{srdefordre}.
Then, by choosing $O_{m}:=O_{m}^{1}\cap O_{m}^{2}$, the
statement dealing with the minimal
order property in Theorem~\ref{srthmA} follows.

Consider a singular trajectory $x(\cdot):=x(\cdot,x_0,u)$ of
$\dot{x}=\sum_{i=1}^m u_i f_i(x)$ admitting
an abnormal extremal $(x(\cdot),\lambda(\cdot))$. Assume that there
exists $J\subset [0,T]\setminus I_{dep}(x(\cdot))$ of
positive measure such that
$G(t)$ is not of rank $m-1$ if $m$ is odd, respectively,
$\widetilde{G}(t)$ is not of rank $m-1$ if $m$ is even.
Following exactly the proofs of Lemmas~3.8 and 3.9 in~\cite{CJT}, the
$m$-tuple  $(f_1,\ldots,f_m)$ belongs to a 
subset of arbitrary codimension in $VF(M)^{m}$ whose complement
contains an open and dense subset. 

We proceed similarly for an argument of the  statement dealing with
the corank one property.  


\subsection{Proofs of Propositions~\ref{propstrictA} 
  and \ref{srpropstrictA}} 
We only treat the control-affine case, the argument for the driftless
control-affine case being identical. We argue by contraposition. 
Consider a nontrivial singular trajectory $x(\cdot):=x(\cdot,x_0,u)$ of
$\dot{x}=f_0(x)+\sum_{i=1}^m u_if_i(x)$, with $(f_0,\ldots,f_m)\in
VF(M)^{m+1}$. Assume that $x(\cdot)$ admits on the one
part a normal extremal lift $(x(\cdot),\lambda^{[n]}(\cdot))$ and on
the other part an abnormal extremal lift
$(x(\cdot),\lambda^{[a]}(\cdot))$. 

Let us introduce some notations. 
For $k\in \N$, let $L=l_1 \cdots l_k$ be a multi-index of
$\{0,\dots,m\}$. The length of $L$ is $|L|=k$ and $f_{L}$ is the 
vector field defined by
$$
f_L:=[[\dots[f_{l_1},f_{l_2}], \dots],f_{l_k}].
$$
A multi-index $L=j l \cdots l$ with $k$ consecutive occurrences 
of the index $l$ is denoted as $L=j l^k$. 

For every multi-index $L$ of $\{0,\ldots,m\}$ and $t\in [0,T]$, set
$$
h_L^{[n]}(t) = \langle\lambda^{[n]}(t),f_L(x(t))\rangle, \
\textrm{and}\ \ 
h_L^{[a]}(t) = \langle\lambda^{[a]}(t),f_L(x(t))\rangle.
$$
After time differentiation, we have on $[0,T]$,
\begin{eqnarray}
\frac{d}{dt} h_{L}^{[n]}(t) 
  = \sum_{l=1}^m u_l(t) h_{Ll}^{[n]}(t)\label{eq:00},\\ 
\frac{d}{dt} h_{L}^{[a]}(t) 
  = \sum_{l=1}^m u_l(t) h_{Ll}^{[a]}(t)\label{eq:11}.
\end{eqnarray}
Recall that, according to the Pontryagin Maximum Principle, there holds
\begin{equation}\label{inconnus}
u(t) = \begin{pmatrix}u_1(t)\\ \vdots \\ u_m(t)\end{pmatrix}
=  U(x(t))^{-1} \begin{pmatrix} h_1^{[n]}(t)
  \\ \vdots \\ h_m^{[n]}(t)\end{pmatrix},
\end{equation}
and, for every $t\in [0,T]$, 
\begin{equation}\label{eqanor}
h_0^{[a]}(t)=\textrm{constant}, \ h_l^{[a]}(t) = 0 ,
\end{equation}
for every $l\in\{1,\ldots,m\}$, and $t\in[0,T]$. Since the trajectory
$x(\cdot)$ is nontrivial, there exists an open interval
$J\subset[0,T]$ and $i\in\{0,\ldots,m\}$ such that
$u_{i}(\cdot)f_{i}(x(\cdot))$ is never vanishing (with the
convention $u_0\equiv 1$). Fix $j\in\{0,\ldots,m\}\setminus\{i\}$.
Differentiating $s$ times (with $s\geq 1$) the relation $h_j^{[a]}(t)
= \mathrm{constant}$ with respect to $t\in J$, one gets, by using
(\ref{eq:00}), (\ref{eq:11}), and (\ref{inconnus}), that
\begin{equation}\label{eq:induc}
0 = \frac{d^s}{dt^s} h_{j}^{[a]}(t) 
  = (u_i(t))^s h^{[a]}_{ji^s}(t) + R_s(t),
\end{equation}
where $R_{s}(t)$ is polynomial in $h_L^{[n]}(t)$ and $h_K^{[a]}(t)$, $\vert
L\vert\leq s$, $\vert K\vert\leq s+1$, with $K$ different from
$ji^s$ and $iji^{s-1}$.
Fix $t\in J$. Since $u_i(t)\neq 0$ and $f_i(x(t))\neq 0$,
we are in a position to apply routine transversality arguments.
It follows that
the $(m+1)$-tuple $(f_0,\ldots,f_m)$ belongs to a closed subset of
$VF(M)^{m+1}$ of arbitrary codimension. Proposition \ref{propstrictA}
follows.


\subsection{Proofs of Propositions~\ref{propstrictB},
  \ref{propstrictB11}, and~\ref{srpropstrictB}} 
We first prove Proposition \ref{propstrictB11} and argue by contraposition.
Consider a nontrivial singular trajectory $x(\cdot):=x(\cdot,x_0,u)$ of
$\dot{x}=f_0(x)+\sum_{i=1}^m u_if_i(x)$. Assume that $x(\cdot)$ admits
on the one part a normal extremal lift
$(x(\cdot),\lambda^{[n]}(\cdot))$ and on the other part an abnormal
extremal lift $(x(\cdot),\lambda^{[a]}(\cdot))$.

 From the Pontryagin Maximum Principle, there holds, for
$l=1,\ldots,m$,
$$ u_l(t) = \sum_{p=1}^m Q^{lp}(x(t))  \beta_p(x(t)),\qquad
\beta_p(x(t)):=h^{[n]}_p(t)-\alpha_p(x(t)),$$
where the $Q^{lp}(x)$ and the $\alpha_p(x)$ are respectively the
coefficients of $U^{-1}(x)$ and of $\alpha(x)$. Note that
the $u_l$'s are smooth functions of the time.

Since the trajectory $x(\cdot)$ is nontrivial, there exists an open
interval $J\subset[0,T]$ such that $\dot x$ is never vanishing on $J$
and one of the two following cases holds.
\medskip

\underline{Case $1$:} $u\equiv 0$ on $J$.

In that case, $\dot x(t)=f_0(x(t))$ for $t\in J$,
and $f_0(x(\cdot))$ is never vanishing on $J$. Moreover, for
$p=1,\ldots,m$, $\beta_p\equiv 0$ on $J$,
i.e. $\alpha_p(x(t))=h^{[n]}_p(t)$ for $t\in J$. By differentiating
the latter relation with respect to the time, we deduce that, for all
$N\geq 0$, $t\in J$ and $p=1,\dots,m$,
$$
L_{f_0}^{N}\alpha_p(x(t))=L_{f_0}^{N}h^{[n]}_p(x(t)),
$$
where $L_{f_0}$ denotes the Lie derivative with respect to the vector
field $f_0$.
Applying routine transversality arguments, it follows that
$\alpha$ belongs to a closed subset of $C^{\infty}(M,\R^m)$
of arbitrary codimension.

\medskip

\underline{Case $2$:} $u$ is never vanishing on $J$.

Using (\ref{s-affine}) and the Lie Algebra Rank Condition, there exist
a multi-index $L$, an index $j_0\in\{0,\ldots,m\}$, and a subinterval of $J$
(still denoted $J$), such that
$$h_L^{[a]}(t)=\mathrm{constant},\ \textrm{and}\ 
h_{Lj_0}^{[a]}(t)\neq 0,$$
for every $t\in J$.
Differentiating $h_L^{[a]}$ on $J$, one gets
\begin{equation}\label{formder}
\begin{split}
0 = \frac{d}{dt}h_L^{[a]}(t) &=
 h_{L0}^{[a]}(t) + \sum_{l=1}^m u_l(t)  h^{[a]}_{Ll}(t) \\
& =  h_{L0}^{[a]}(t) + \sum_{1\leq l\leq p\leq m}
 c_{lp}(t)  Q^{lp}(x(t)) ,
\end{split}
\end{equation}
where $c_{ll}(t) :=  \beta_l(t) h^{[a]}_{Ll}(t) $, and
$ c_{lp}(t) := \beta_p(t) h^{[a]}_{Ll}(t) +  \beta_l(t)
h^{[a]}_{Lp}(t)$ if $l< p$.

\begin{lemma}
Up to reducing the interval $J$, there exist indices $j$ and $l$ in
$\{1,\ldots,m\}$ such that $c_{lj}(t)$ or $c_{jl}(t)$ is never
vanishing on $J$.
\end{lemma}

\begin{proof}
If $j_0=0$, then $h^{[a]}_{L0}(t)\neq 0$, and it follows from
(\ref{formder}) that there exist $l,j\in\{1,\ldots,m\}$ such that
$c_{lj}(t)\neq 0$. Otherwise, take $j:=j_0$. In that case, one of the
$\beta_p$'s
does not vanish on $J$ since $u$ is not zero. First, assume that 
$\beta_{j}(t)$ is not identically equal to zero on $J$; then, up to
reducing $J$, 
$c_{jj}(t)$ is never vanishing on $J$. Otherwise, there exists $l\neq
j$ such that, up to reducing $J$, $\beta_l$ is never vanishing on $J$
and thus similarly for $c_{lj}$ (or $c_{jl}$).
\end{proof}

For $t\in J$, let $F_t\in VF(M)$ be the vector field defined by
$$F_t(x) := f_0(x)+ \sum_{i=1}^m u_i(t) f_i(x).$$
Notice that $F_t(x(t))=\dot{x}(t)\neq 0$.
For all $N\geq 0$ and $t\in J$, we get, by taking the $(N+1)$-th time
derivative of $h_L^{[a]}$ on $J$, 
$$
0 = \frac{d^{N+1}}{dt^{N+1}}h_L^{[a]}(t) =
 c_{lj}(t) L^N_{F_t} Q^{jl}(x(t))
+ R_N(t),
$$
where $R_N(t)$ is a linear combination of $L^s_{F_t}
Q^{pi}(x(t))$, with $s\leq N$, $p\leq i$ in $\{1,\dots,m\}$ and
$s<N$ if  $(p,i)=(j,l)$, and of
$L^s_{f_r}Q^{pi}(x(t))$, with $s<N$, $p\leq i$ in $\{1,\ldots,m\}$,
and $r\in\{0,\ldots,m\}$.
Applying routine transversality arguments, 
it follows that $(U,\alpha)$ belongs to a closed subset of 
$\mathcal{S}^+_m(M)\times C^{\infty}(M,\R^m)$
 of arbitrary codimension. Proposition
\ref{propstrictB11} is proved.

To show Propositions~\ref{propstrictB} and~\ref{srpropstrictB}, simply
notice that the 
argument of Case $2$ with $\alpha=0$ applies with suitable
modifications.  


\end{document}